\documentclass{article}

\begin{document}

\title{Integral bounds on curvature and Gromov-Hausdorff limits}
\author{X-X. Chen and S. K. Donaldson}
\maketitle


\newtheorem{thm}{Theorem}
\newtheorem{prop}{Proposition}
\newtheorem{lem}{Lemma}
\newtheorem{defn}{Definition}
\newcommand{\reg}{{\rm reg}}
\newcommand{\Ric}{{\rm Ric}}
\newcommand{\bR}{{\bf R}}
\newcommand{\bO}{{\bf O}}
\newcommand{\Vol}{{\rm Vol}}
\newcommand{\Riem}{{\rm Riem}}
\newcommand{\cA}{{\cal A}}
\newcommand{\uf}{\underline{f}}
\newcommand{\tuf}{\tilde{\underline{f}}}
\newcommand{\ugamma}{\underline{ \gamma}}
\section{Introduction}

In this note we develop a new approach to certain results of Cheeger, Colding and Tian involving the curvature of Riemannian manifolds close, in the Gromov-Hausdorff sense, to a singular limit. In one aspect our results go a little beyond those in the literature but the main interest, for us, is that the arguments are substantially different  and may be more amenable to certain generalisations.

Let $n,k$ be positive integers with $n\geq k>2$. Let $\Gamma\subset SO(k)$ be a non-trivial finite group which acts freely on the unit sphere $S^{k-1}$. Write $M_{n,\Gamma}$ for the singular space $\left( \bR^{k}/\Gamma \right)\times R^{n-k}$. Let $B_{n,\Gamma}$ be the unit ball centred at the equivalence class of $(0,0)$. In other words, $B_{n,\Gamma}$ is the quotient of the unit ball in $\bR^{n}$ by the action of $\Gamma$, embedded in $SO(n)$ in the obvious way.

We consider a complete Riemannian manifold $(X^{n},g)$, with bounded Ricci curvature: $\vert \Ric\vert\leq (n-1)$. (In practice $\vert \Ric \vert$ will usually be very small, by re-scaling.) We assume that we have a non-collapsing condition:
\begin{equation}  \Vol (B(x,r))\geq C^{-1} r^{n} \end{equation}
for all balls in $X$, when $r$ is less than the diameter of $X$. Then we have
\begin{thm}
There are $\delta=\delta(n,C)>0$ and $\zeta=\zeta(n,C)>0$ such that if the  Gromov-Hausdorff distance from a unit ball $B(x_{0},1)$ in $X$  to $B_{n,\Gamma}$  is less than $\delta$ then 
$$  \int_{B(x_{0},1)} \vert \Riem \vert^{k/2} \geq \zeta. $$

\end{thm}

Next let $0<\beta<1$ and let $\bR^{2}_{\beta}$ be the standard cone with cone angle $2\pi \beta$. Thus in polar co-ordinates the metric is
$dr^{2}+ \beta^{2} r^{2} d\theta^{2}$. For $n\geq 2$ let $N_{n,\beta}$ be the product $\bR^{2}_{\beta}\times \bR^{n-2}$ and let $B_{n,\beta}$ be the unit ball centred at the origin, in the obvious sense.
\begin{thm}
There are $\delta=\delta(n,C,\beta)>0$ and $\zeta=\zeta(n,C,\beta)>0$ such that if the  Gromov-Hausdorff distance from a unit ball $B(x_{0},1)$ in $X$  to $B_{n,\beta}$  is less than $\delta$ then 
$$  \int_{B(x_{0},1)} \vert \Riem \vert \geq \zeta. $$
\end{thm}

In the case of  K\"ahler metrics, Theorem 1 was proved by Cheeger, Colding and Tian in \cite{kn:CCT}.  They also establish the result in the real case for many groups $\Gamma$. See also results of Cheeger and Tian in \cite{kn:CT}. In the co-dimension 2 case (Theorem 2), Cheeger proved in \cite{kn:C} a significantly stronger result than that stated here,  involving only a lower bound on the Ricci curvature. In fact a conjecture of Anderson \cite{kn:A2} would, if true, mean that, in Theorem 2, if $\delta$  is sufficiently small there is no such ball in any $X$.  It seems possible to us that the approach here could be extended to prove the analogues of Theorem 1 and 2 for metrics with Ricci curvature bounded below but this certainly involves significant new problems.

Our proof of Theorem 1  will use the co-dimension 2 result, Theorem 2. We also sketch a proof of Theorem 2, which is similar but  involves some  extra complications. So we will discuss Theorem 1 first,  in Sections 2, 3 and 4 below, and then discuss Theorem 2 in Section 5.  Of course we can also quote what we need from \cite{kn:C}.

The central notion in our proof is that of a point $x$ in $X$ around which the \lq\lq energy is small at all scales''. For any $r$ we set
$$   E(x,r)= r^{2-n} \int_{B(x,r)} \vert \Riem\vert. $$
  To explain the point of the definition: given a ball $B(x,r)$ write $B(x,r)^{\sharp}$ for the same ball with the metric (i.e lengths) scaled by a factor $r^{-1}$. So $B(x,r)^{\sharp}$ is a unit ball in the new metric. Then $E(x,r)$ is the $L^{1}$ norm of the curvature of $B(x,r)^{\sharp}$. Set $\overline{E}(x)=\max_{r\leq 2} E(x,r)$. (The restriction to $r\leq 2$ is just because we are concerned with the local picture.) For $\epsilon>0$ let $\cA_{\epsilon}\subset X$ be the set of points $x$ where $\overline{E}(x)\leq \epsilon$. Thus
$\cA_{\epsilon}$ is the set of points where \lq\lq the energy is less than $\epsilon$ at all scales''(or, more precisely, at all small scales). For the proof of Theorem 1 we could work with a similar definition using $L^{k/2}$ norms. A disadvantage of that is that {\it a posteriori} the set corresponding to $\cA_{\epsilon}$ would  be empty  for small enough $\epsilon$ and the statements we prove about it would be vacuous. That is one reason why we prefer to use $L^{1}$. (A similar issue arises in the case of Theorem 2, see the discussion at the end of  Section 5 below).    

With this background we can attempt to explain the central idea of our proof in very vague informal terms, but perhaps sufficient that an expert could reconstruct   the arguments. Suppose $\delta$ is very small and we start at distance roughly $1$ from the \lq\lq approximately singular set''. Then the metric appears very close to a flat cone with a codimension $k$  singularity.  As we travel towards the  the approximately singular set, viewing at smaller and smaller scales,  the apparent singularity must eventually resolve into a smooth metric and at this scale we must see some curvature, measured in $L^{1}$ norm on a ball of the appropriate scale. Conversely at a point in $\cA_{\epsilon}$ where the energy is small at all scales we still see the same apparent singularity--or a flat space--at every scale. That is, staying in $\cA_{\epsilon}$ we can never reach the apparent singularity. This implies that the complement of $A_{\epsilon}$ must be \lq\lq large'' in a  certain topological sense and that yields the desired lower bound by a standard covering argument.

We should emphasise that while our proofs, in so much as they are new, are elementary and self-contained they depend strongly on deep foundational results  in this area due to Anderson, Cheeger, Colding and Tian.
\section{Set up and strategy}

We put some notation and conventions in place for our later arguments. First to simplify exposition slightly we will first give the proof of Theorem 1 in the case of an Einstein manifold $X$, so $\Ric=\lambda g$ with $\vert \lambda\vert \leq 1$. Then in Section 4.2 we explain the arguments needed to extend to bounded Ricci curvature.
We suppose throughout Sections 3 and 4 that we have some fixed $\Gamma\subset SO(k)$. (For fixed $n$ and $C$ there are only finitely many possible choices of $\Gamma$, up to conjugacy.) Note first that the complete manifold $X$---far away from the point $x_{0}$---- plays no real role since all our arguments will be local. Similarly there is no loss in supposing that a substantially larger ball $B(x_{0}, K)$ (with $K=10$, say) is also Gromov-Hausdorff close to the $K$-ball  in the model. This is because we could always restrict attention to a smaller ball and re-scale.  With this understood, we do not need to distinguish between the distance function in $B(x_{0},1)$ and that induced from the metric in $X$.
We write a point in $M_{n,\Gamma}$ as a pair $(\xi,\eta)$ with $\xi\in \bR^{k}/\Gamma$ and $\eta\in \bR^{n-k}$. We write $\vert \xi\vert$ with the obvious meaning. For $s>0$, let $H_{s}$ be the set 
$$   H_{s}= \{ (\xi,\eta): \vert \xi \vert^{2}+\vert \eta\vert^{2}\leq 1-s \ , \ \vert \xi \vert\geq s\}. $$

Let $\Omega\subset H_{s}$ be the region where $\vert \eta\vert \leq (1/10)+ \vert \xi\vert/5$. Thus $\Omega$ meets the outer boundary $\vert \xi\vert^{2}+\vert\eta\vert^{2}=1-s$ in an \lq\lq equatorial region'' $E$. The crucial property is that $\xi$ is never zero on $E$.

Now recall that to say that the  Gromov-Hausdorff distance between $B(x_{0},1)$ and $B_{n,\Gamma}$ is less than $\delta$ is to say that   there is a metric on the disjoint union which extends the given metrics on the two subsets and such that the $\delta$-neighbourhood of either subset is the whole of the union. In our situation we can use a much more concrete notion. Define an
{\it $s$-chart} at $x_{0}$  to be  a map $\chi_{s}: H_{s}\rightarrow B(x_{0},1)$ with the following properties
\begin{itemize}
\item  $\chi_{s}$ is a diffeomorphism to its image.
\item The pull-back of the metric $g$ differs  in $C^{2}$ norm  from the given flat metric on $H_{s}$ by at most $s$.
\item $ B(x_{0},1)$ is contained in the $2s$ neighbourhood of the image of $\chi_{s}$.
\end{itemize}
Then  under our hypotheses it is a fact that for any $s$ there is an $s$-chart if $\delta$ is small enough. This depends on a great deal of deep theory, notably Anderson's results on the volume ratio \cite{kn:A} and Colding's result on volume convergence \cite{kn:Co}. We refer to  \cite{kn:C2}, \cite{kn:C3}  for  accounts of the whole theory.  Of course an $s$-chart is not unique but it is essentially unique up to the isometries of $B_{n,\Gamma}$ and a small arbitrary error. Recall that for the time being we are considering Einstein metrics. This means that in the second item we could replace $C^{2}$ by $C^{r}$ for any $r$ using elliptic regularity of the Einstein equations in local harmonic co-ordinates.

Given an $s$-chart, let $\Omega'\subset B(x_{0}, 1)$ be the union of
$\chi_{s}(\Omega)$ and the ball $B(x_{0},1/10)$. This has an \lq\lq outer boundary''given by $\chi_{s}(E)$. A moments thought will show the reader that, in proving our main theorem, there is no loss of generality in assuming, for any given $\epsilon$, that $\chi_{s}(E)$ lies in $\cA_{\epsilon}$. (Here again we use the freedom to restrict to a smaller ball in the original problem.) 
 With all this preparation we can state the central result in our proof.
 \begin{prop}
 There is an $s_{0}$ and for all $s\leq s_{0}$ an $\epsilon(s)$ such that if an $s$-chart exists and $\epsilon<\epsilon(s)$ then there is a continuous retraction $R: \Omega'\cap \cA_{\epsilon}\rightarrow \chi_{s}(E)$ which restricts
to the identity map on $\chi_{s}(E)$.
\end{prop}

This is, in precise language, what we mean by saying that \lq\lq  staying in  $\cA_{\epsilon}$ we can  never reach the singularity''.

In the remainder of this section we will explain the proof that Proposition 1  implies Theorem 1. This involves the construction of \lq\lq slices''. What we need is quite standard and elementary. (There is some similarity here with the approach of Cheeger, Colding and Tian but they need a much more sophisticated slicing, controlling higher derivatives.) Recall that we are supposing   that the metric is Gromov-Hausdorff close, say distance $\delta$,  to the flat model on the much larger ball $B(x_{0}, K)$.
For $i=1,\dots, n-k$ let $p_{i}$ be the point in $\bR^{n-k}$ with ith. coordinate $K$ and all others $0$. Choose a point $q_{i}\in B(x_{0},K)$ which is a distance less than $\delta$ from $(0,p_{i})$. Let $f_{i}$ be the function on $B(x_{0},1)$
$$   f_{i}(x)= d(x_{0}, q_{i})- d(x,q_{i})$$ and let $\uf:B_{x_{0}}(1)\rightarrow \bR^{n-k}$ be the map with components $f_{i}$. Clearly for the corresponding construction in the flat model this map approaches the projection from $\left(\bR^{k}/\Gamma\right)\times \bR^{n-k}$ to $\bR^{n-k}$ as $K$ tends to infinity. It follows that if $K$ is sufficiently large (how large being something one can determine by elementary geometry in the flat model)  then we can find a $\alpha>0$ such that if we have an s-chart with sufficiently small $s$ then for any $\eta$ with $\vert \eta\vert \leq \alpha$ the fibre $\uf^{-1}(\eta)$ is contained in  $\Omega'$. 

The function $\uf$ is Lipschitz, with Lipschitz constant $1$. It is then standard that we can find a smooth map $\tuf$, arbitrarily close to $f$ in $C^{0}$ with derivative bounded in operator norm by $2$ say (or any number bigger than $1$). It is also clear that we can suppose that   near the boundary region $\chi_{s}(E)$ the composite $\tuf \circ \chi_{s}$ is exactly equal to the projection map from $\left( \bR^{k}/\Gamma\right)\times
\bR^{n-k}$ to $\bR^{n-k}$ and that fibres $\tuf^{-1}(\eta)$ are contained in $\Omega'$, for $\vert \eta\vert\leq \alpha$. 

Now let $\pi:\chi_{s}(E)\rightarrow \bR^{k}/\Gamma$ be the composite of $\chi_{s}^{-1}$ with the map $(\xi,\eta)\mapsto \xi/\vert \xi\vert$. For almost all $\eta$ with $\vert \eta\vert \leq \alpha$ the fibre $\tuf^{-1}(\eta)$ is a compact $k$ manifold with boundary and $\pi$ yields a diffeomorphism from the boundary to $S^{k-1}/\Gamma$. This means that $\tuf^{-1}(\eta)$ cannot be contained in $\cA_{\epsilon}$, for then the map $R\circ \pi$ would lead to a retraction of $\tuf^{-1}(\eta)$ onto its boundary. To sum up, we have established the following.
\begin{prop}
For suitable choice of $\delta,\epsilon$ there is an $\alpha>0$ and a smooth map $\tuf$ from $B_{x_{0}}(1)$  to $\bR^{n-k}$ with derivative bounded by $2$ and which maps the complement of $\cA_{\epsilon}$ onto the $\alpha$-ball in $\bR^{n-k}$. 
\end{prop}
 
 Now the proof is completed by a standard Vitali covering argument(as in  \cite{kn:S}, Section 1.6 for example).
Call a ball $B(x,r)$ with $E(x,r)>\epsilon$ a \lq\lq high energy ball''.
We successively choose high energy balls $B_{1}, B_{2},\dots $ such that at stage $i$ the ball $B_{i}=B(x_{i}, r_{i})$ is disjoint from $B_{1}\dots, B_{i-1}$ and has maximal radius among all such possibilities. (If there are no such possibilities we stop.)  Directly from the definition, any point of the complement of $\cA_{\epsilon}$ is the centre  of some high energy ball. Then the selection scheme implies that the complement of $\cA_{\epsilon}$ is contained in the union of the closures of the twice-sized balls $B(x_{i}, 2r_{i})$. 
Thus $\tuf$ maps $\bigcup_{i}\overline{ B(x_{i}, 2r_{i})}$ onto the $\alpha$- ball in $\bR^{n-k}$. Write $\omega_{n-k}$ for the volume of the unit ball in $\bR^{n-k}$. Then the derivative bound on $\tuf$ implies that the volume of $\tuf B(x_{i}, 2r_{i})$ is at most $\omega_{n-k} (2 r_{i})^{n-k}$. Since the volume of the $\alpha$ ball is $\omega_{n-k} \alpha^{n-k}$, we get
\begin{equation}   2^{n-k} \sum r_{i}^{n-k} \geq \alpha^{n-k} . \end{equation}
On the other hand since the balls $B_{i}$ are disjoint we have
$$  \sum \int_{B_{i}}\vert \Riem\vert^{k/2} \leq \int_{B_{x_{0}}(1)} \vert \Riem \vert^{k/2}. $$
By the definition of a high energy ball
$$  \int_{B_{i}} \vert \Riem \vert \geq \epsilon r_{i}^{n-2}. $$
The bound on the Ricci curvature implies an upper bound on the volume of each of these balls, say $\Vol(B_{i})\leq c r_{i}^{n}$. Then H\"older's inequality gives
$$   \int_{B_{i}} \vert \Riem \vert^{k/2} \geq c^{1-k/2}\epsilon^{k/2}  r_{i}^{n-k}. $$

Using (2), we conclude that
$$  \int_{B(x_{0},1)} \vert \Riem \vert^{k/2}\geq c^{1-k/2}\epsilon \alpha^{n-k} 2^{k-n}. $$
\section{The main proof}
\subsection{Strategy}
The purpose of this Section is to prove Proposition 1. Our approach is to construct a vector field $v$ on $B(x_{0}, 1)$ such that the integral curve of $v$ starting at any point $x$ in $\cA_{\epsilon}\cap \Omega'$ hits $\chi_{s}(E)$  after some finite time and then define $R(x)$ to be the hitting point. In terms of our chosen  $s$-chart $\chi_{s}$ on $B(x_{0},1)$ it is obvious how one might do this if one only considers points $x$ which also lie in the image $\chi_{s}(H_{s})$.
That is, in the flat model,  define $\partial_{r}$ to be the unit radial vector field , given by the gradient of the function $\vert \xi\vert$. Then  we could simply use the push forward of $\partial_{ r}$ under $\chi_{s}$. The whole problem is to show that we can extend this definition---or something like it---outside $\chi_{s}(H_{s})$.
This divides into two parts. One is to show that for a point in $\cA_{\epsilon}$ we see essentially the same picture at any scale. The ingredients for that are developed in this subsection 3.1, except for the proof of Theorem 3 below, which we postpone to Section 4.  The other part is to devise some method of defining the vector field $v$ which, roughly speaking, looks like $\partial_{r}$ at all scales. The ingredients for that are developed in 3.2. In 3.3 we bring these ideas together to prove Proposition 1  

We will say that a space $Z$ is an  {\it $L^{1}$-flat limit ball} if \begin{itemize}
\item Z is the Gromov-Hausdorff limit of unit balls $B(p_{i}, 1)\subset X_{i}$ where $X_{i}$ are Einstein $n$-manifolds satisfying the conditions we considered in Section 1 (with a fixed $C$ in (1));
\item For some $r>1$ 
$$   \int_{B(p_{i}, r)} \vert \Riem \vert \rightarrow 0 , $$
as $i\rightarrow \infty$.
\end{itemize}
 Thus we know that an open dense subset $Z_{\reg}$ of $Z$ has
a flat Riemannian metric. The definition on a $s$-chart extends in the obvious way to this situation---we just require that $\chi_{s}$ maps in to the smooth subset. We have the following rigidity result, proved in Section 4.  

\begin{thm}
Suppose $B$ is an $L^{1}$-flat limit ball which admits an   $s$-chart. If $s$ is   sufficiently small then $B$ is isometric to a ball in $M_{\Gamma}$.
\end{thm}

Let $\Sigma$ be the set of points of the form $(\xi,0)$ in $H_{s}\subset M_{\Gamma,n}$. (Thus, roughly speaking, $\Sigma$ is a transversal to the singular set.) Recall that we are using the notation $B^{\sharp}$ for  rescaled balls and we will write the rescaled metric as $d^{\sharp}$.
\begin{prop}
Suppose $B(p,1)$ is a ball in the manifold $X$ which admits an $s$-chart and let $q$ be a point with $d(q,p)\leq 1/10$. For $s$ sufficiently small there is an $\epsilon(s)>0$ such that if the $E(p,1)<\epsilon(s)$then  there is a point $p'\in B(p,1)$ with $d(p',q)<1/2$ and such that the rescaled ball $B(p',1/2)^{\sharp}$ admits an $s$ chart $\chi_{s}$. This $s$-chart can be chosen such that  if $d^{\sharp}(p',q)\geq 1/10$ then $q$ lies in the image $\chi_{s}(\Sigma)$.
\end{prop}

It is easy to see here that $d^{\sharp}(p',q)$ cannot be more than $1/5+O(s)$.
This proposition expresses, in precise language, what we mean by saying that
\lq\lq at a point in $\cA_{\epsilon}$ \dots we still see the same apparent singularity, or a flat space, at every scale''.  

To prove this proposition we argue  by contradiction. We can suppose that $s$ is less than the value given by  Theorem 3. If the statement fails we get a sequence of such balls with $L^{1}$ norm of the curvature tending to zero.  We get a flat limit ball 
which admits an $s$-chart and hence, by the theorem, is itself isometric to a ball in $M_{\Gamma,n}$. Thus the half-sized ball in the limit is isometric to a ball in $M_{\Gamma,n}$ and this gives a contradiction. The condition that $q$ lies in the image of $\chi(\Sigma)$, if $d^{\sharp}(q,p')>1/10$ is achieved by making a small translation in the $\bR^{n-k}$ factor in $M_{\Gamma,n}$.

\subsection{The lasso function}
               Given any ball $B\subset B(x_{0},1)$ such that the normalised ball $B^{\sharp}$ admits an $s$ chart, for some small $s$, it is clear, as we said in the previous subsection,  how we could define an appropriate vector field on the image of the s-chart in $B$. The problem is that when $B$ is a tiny ball not contained in the image of the $s$-chart in $B(x_{0},1)$ (which is the case of interest) there is no {\it a priori} notion of compatability between the $s$ chart in $B$ and the s-chart in $B(x_{0},1)$. One way around this might be to build up a global vector field by a patching construction but instead we prefer to concoct a global definition and show that it behaves in the right way at every scale. This definition requires a digression.
\begin{defn}
Let $V$ be a complete Riemannian manifold and fix $\lambda>0$. For $p\in V$ a lasso based at $p$ is a pair $(\gamma_{0}, \gamma_{1})$ where $\gamma_{i}:[0,l_{i}]\rightarrow V$ are geodesics segments parametrised by arc length such that 
\begin{itemize}
\item $l_{1}>0$;
\item $\gamma_{0}(0)=p$;
\item $\gamma_{0}(l_{0})=\gamma_{1}(0)=\gamma_{1}(l_{1})$.
\end{itemize}
We call $\gamma_{0}$ the {\it knot} of the lasso and define the $\lambda$-length of the lasso to be $l_{0}+ \lambda l_{1}$. If $l_{0}>0$ we will say the lasso is a {\it genuine lasso}.
\end{defn}

If $V$ is compact (say) it is clear that the set of lassos based at $p$ is not empty and we thus have a function $\rho_{\lambda}(p)$ defined by the infimum of the $\lambda$-length. It follows immediately from the definition that $\rho_{\lambda}$ is Lipschitz, with Lipschitz constant $1$. It is also clear that the infimum  is realised by at least one {\it minimising lasso}. For a genuine lasso, the minimising condition implies (by the usual first variation formula about geodesics) that
\begin{equation}  \gamma_{0}'(l_{0}) = \lambda\left( \gamma_{1}'(0)- \gamma_{1}'(l_{1})\right). \end{equation}
Thus the angle between the two tangent vectors of $\gamma_{1}$ at the knot  is determined by $\lambda$. If $\lambda<1/2$ then this condition (3) can never be realised, so a minimising lasso is never genuine. In this case we have
$\rho_{\lambda}(p)= 2 \lambda I(p)$, where $I(p)$ is the injectivity radius at the point $p$. In fact in our application we could use this injectivity radius function but the variant $\rho_{\lambda}$ (for a suitable large value of $\lambda$) will make the argument simpler and more transparent. For any minimising lasso $\ugamma=(\gamma_{0},\gamma_{1})$ based at $p\in V$ we define the {\it tension vector} $\tau_{\ugamma}\in TV_{p}$ to be $-\gamma_{0}'(0)$ if the lasso is genuine and $\lambda( \gamma_{1}'(l_{1})- \gamma_{1}'(0))$ if not. (In our situation we will only really be concerned with the case of genuine lassos.)

The function $\rho_{\lambda}$ will not necessarily be smooth. We make a further digression to review a useful general principle by which we can get around this difficulty. Suppose $P,Q$ are manifolds and $F$ is a smooth function defined on an open set $U\subset P\times Q$. For simplicity we assume that $Q$ is compact. Let $\pi_{Q}:U\rightarrow Q$ be the restriction of the projection map and for each $q\in Q$ let $U_{q}=\pi_{Q}^{-1}(q)$. Let
$f(q)$ be the infimum of $F$ on $U_{q}$. Suppose that this is finite and that there is a compact subset $K_{q}\subset U_{q}$ and $\delta_{q}>0$ such that $ F\geq f(q)+\delta_{q}$ on $U_{q}\setminus K_{q}$. Thus the infimum
$f(q)$ is attained on a compact set $J_{q}\subset K_{q}$. It is clear that the function $f$ on $Q$ is continuous, but in general it may not be smooth. For each point $(p,q)$ in $J_{q}$ let $(D_{Q}F)(p,q)$ denote the partial derivative of $F$ at  $(p,q)$ in the $Q$ variable. This can be regarded as an element of $T^{*}Q_{q}$.
 So we have a map, say $\iota_{q}:J_{q}\rightarrow T^{*}Q_{q}$. Take the convex hull of the image in the vector space $T^{*}Q_{q}$. Then as $q$ varies over $Q$ the union of these yields a  subset ${\cal D}\subset T^{*} Q$. In the next proposition we will regard the derivative of a function on $Q$ as a subset of $T^{*} Q$.
 
 \begin{prop}
 In this  situation $f$ can be approximated arbitrarily closely in $C^{0}$ by a smooth function whose derivative lies in an arbitrarily small neighbourhood of ${\cal D}$. 
\end{prop}
To see this, we can first approximate $F$ by a function $\tilde{F}$ whose restriction to slices $U_{q}$ has nondegenerate critical points, for  $q$ outside a set $\Delta\subset Q$ of codimension at least $1$.  We can further suppose that the minimum is unique, for $q$ outside $\Delta$. The minimiser lies arbitrarily close to $J_{q}$. Define $\tilde{f}$ on $Q$  using $\tilde{F}$ in place of $F$.  Outside $\Delta$ it is clear that $\tilde{f}$ is smooth and by elementary calculus its derivative is given by the partial derivative of $\tilde{F}$ at the minimiser. Then smooth $\tilde{f}$ in the standard way using a suitable family of integral operators 
$$  \tilde{f}_{\epsilon}= \int_{Q} \beta_{\epsilon}(q,q') \tilde{f}(q') dq'. $$
It is straightforward to verify that these approximations have the desired properties.

\

 Now return to the lasso function $\rho_{\lambda}$ on $V$. For each $p\in V$ we take  the convex hull of the tension vectors
$\tau_{\ugamma}$ as $\ugamma$ runs over all the mimimising lassos based at $p$ and let ${\cal D}\subset TV $ be the union of these, much  as above. Then we have

 \begin{prop}
 The function $\rho_{\lambda}$ on $V$ can be approximated arbitrarily closely in $C^{0}$ by a smooth function whose gradient lies in an arbitrarily small neighborhood  of ${\cal D}$.
 
\end{prop}

 To see this we write the length of $\gamma_{i}$ as the square root of the \lq\lq energy''. There is a minor complication because the square root is not smooth when $\gamma_{0}$ has length zero. (In fact this case will not enter in our application below, so we will be rather sketchy.) We choose a family of functions    $W_{\eta}(x)$ (for $\eta>0$), with $W_{\eta}(x)$ a smooth function of $x^{2}$ in a small neighbourhood of $x=0$ and $W_{\eta}(x)=x$ for $x\geq \eta$. Now perturb the length function
to $$ W_{\eta}(l_{0})+ \lambda l_{1}$$ and consider lassos which minimise this perturbed functions. The appropriate matching condition, generalising (3),  is
$$   W_{\eta}'(l_{0})
\gamma'_{0}(l_{0}) = \lambda(\gamma'_{1}(0)- \gamma'_{1}(l_{1}))$$
where the left hand side is interpreted as zero if $\gamma_{0}$ is a constant map.
With the same convention, we define the tension vector to be $\tau_{\ugamma,\eta}= W_{\eta}'(l_{0})\gamma'_{0}(0)$.  We want to prove the obvious analogue of Proposition 5  for the perturbed problem. Once this is done the proof of Proposition 5 follows easily by letting $\eta$ tend to zero.
 Just as in the usual Morse Theory of geodesics based on the energy functional (\cite{kn:M}), we can construct a finite-dimensional space of \lq\lq piecewise geodesic  lassos'' which can be identified with an open set in the product $V^{N}$ of a large number of copies of $V$. Here we parametrise curves by the fixed interval $[0,1]$. The minimisers for the perturbed problem are the minimisers of a function $\tau(\sqrt E_{1}) + W_{\eta}(\sqrt{E_{0}})$ where $E_{0}, E_{1}$ are the energies of the two piecewise-geodesic curves. We fit into the framework of Proposition 4, with $P=V$ and $Q=V^{N-1}$, and the result follows.  
\subsection{The main argument}
In this subsection we bring together the strands developed above. We begin by considering the flat model space $M_{\Gamma,n}$. Consider  a geodesic loop in the smooth part of  $M_{\Gamma,n}$. This lifts to a line segment in $\bR^{n}$  joining two points in the same $\Gamma$ orbit. It is clear that the angle between the two tangent vectors at the base point cannot be too small. So if we choose the parameter $\lambda$ large enough (how large depending solely on $\Gamma$) then the matching condition (3) can never be satisfied. Next let $\nu$ be a point in $\bR^{k}/\Gamma$ with $\vert \nu\vert=1$.
It is clear that there is some $c>0$ such that, for all such $\nu$ there are distinct geodesics segments in the smooth part of $\bR^{k}/\Gamma$ emanating  from $\nu$, of lengths less than $c$, with the same endpoints and whose distance from the origin is greater than $c^{-1}$. 

Suppose that $B$ is a ball in $X$ centred at $p$, such that the normalised ball $B^{\sharp}$ admits an $s$-chart $\chi_{s}$.Recall that we write $d^{\sharp}$ for the normalised metric. Simplifying notation, we will  write $\rho$ for  $\rho_{\lambda}$ on $X$ and $\rho^{\sharp}$ for the normalised version in $B^{\sharp}$.  We will now ignore the fact that $\rho$ is not smooth, since for our purposes below we can always pass to a smooth approximation. 

Let $x=\chi_{s}(\xi,\eta)$ be a point in the image of $\chi_{s}$. Define $$  \varpi(x)= \chi_{s}( cs \vert \xi\vert^{-1} \xi, \eta). $$
Then the preceding discussion implies that the injectivity radius at $\varpi(x)$ (for the normalised metric) is not more than $c^{2}s$. Since the distance from $d^{\sharp}(x,\varpi(x))$ is at most $\vert \xi \vert +O(s)$ we see that there are lassos based at $x$ of $\lambda$-length (in the normalised metric) less than $\vert \xi \vert + C s$, for some fixed $C$. Thus $\rho^{\sharp}_{\lambda}(x)\leq \vert \xi \vert + Cs $. Now fix a constant $\kappa_{1}$ slightly less than $1/2$ and suppose that $d^{\sharp}(p,x)\leq \kappa_{1}$. Then it follows that, if $s$ is sufficiently small, any lasso  based at $x$ of $\lambda$-length less than $\vert \xi\vert + Cs$ must lie in $B^{\sharp}$. On the other hand, by our choice of $\lambda$ the lasso cannot lie entirely in the image of $\chi_{s}$. It is then clear that $\rho^{\sharp}(x)= \vert \xi\vert+ O(s)$. If we consider points with  $\vert \xi \vert$  bounded below by some fixed number then it is clear that (for small $s$)  any minimising lasso based at $x$ is genuine and has a knot a (normalised) distance $O(s^{1/2})$ from $\varpi(x)$. Then  Proposition 5 shows that the gradient of $\rho^{\sharp}$ at $x$ differs from $\partial_{r}$ (interpreted via the $s$-chart) by $O(s^{1/2})$. In particular what we shall need is that for any $\kappa_{0}>0, \kappa_{3}>1$ we can suppose (by making $s$ small) that at points $x$ with $d^{\sharp}(x,p)\leq \kappa_{1}, \rho(x)\geq \kappa_{0}$ we have
\begin{equation} \vert {\rm grad} \rho^{\sharp} \vert \geq \kappa_{3}^{-1}\end{equation}
We fix $\kappa_{0}$ slightly smaller than $1/10$.

Now fix $\kappa_{3}$ slightly bigger than $1$ and choose $\kappa_{2}>0 $ so that \begin{equation}\frac{1}{5}+ 2\kappa_{3}\kappa_{2}<\kappa_{1}\end{equation}
Thus $\kappa_{2}$ is approximately $3/20$ and in particular $\kappa_{2}>\kappa_{0}$.
Let $q$ be a point in $ B^{\sharp}$ with $d^{\sharp}(p,q)< \kappa_{1}-2\kappa_{2}\kappa_{3}$
 We define a set
\begin{equation}  S_{q}=\{ x\in B^{\sharp}: d^{\sharp}(x,q)< \kappa_{3} \rho^{\sharp}(x) \}.\end{equation}
Then we have 
\begin{lem}
Let $x_{0}\in S_{q}$ and $\kappa_{0}\leq \rho^{\sharp}(x_{0})\leq  \kappa_{2}$. Then there is a time interval $[0,T]$ so that gradient flow
$x_{t}$ of ${\rm grad} \rho^{\sharp}$ starting with $x_{0}$ is defined and lies in $S_{q}\subset B^{\sharp}$ for $t\leq T$ and such that $\rho^{\sharp}(x_{T})= 2 \kappa_{2}$.
\end{lem}
To see this, observe that the flow is defined and lies in $S_{q}$ for at least a short time. While the flow is in $S_{q}$ and while $\rho^{\sharp}(x_{t})\leq 2\kappa_{2}$ we have $d^{\sharp}(x,p)\leq d^{\sharp}(x,q)+ d^{\sharp}(q,p)\leq 1/5 + 2\kappa_{3}\kappa_{2}<\kappa_{1}$ by our choices. Since $\rho^{\sharp}(x_{t})\geq \rho^{\sharp}(x_{0})\geq \kappa_{0}$ the gradient bound (4) holds. This implies that, over this time interval, the function $d^{\sharp}(x_{t},q)-\kappa_{3}\rho^{\sharp}(x_{t})$ is decreasing, so we can never reach the boundary of $S_{q}$ and the result follows.

Now  we can give the main argument. Write $B_{0}$ for our original unit ball $B(p_{0},1)$ which admits an $s$-chart for suitably small $s$. Suppose  $q_{0}$ is a point in $B_{0}$ with  $d(p,q_{0})\leq 1/10$ and that $q_{0}$ is in $\cA_{\epsilon}$ where $\epsilon=\epsilon(s)$ as determined by Proposition 3. Applying that Proposition we get another centre $p_{1}$ such that the rescaled ball $B(p_{1}, 1/2)^{\sharp}$ admits an $s$-chart with the same value of $s$. If $d^{\sharp}(p_{1}, q_{0})>1/10$ we stop, otherwise we can repeat the discussion to get a new centre $p_{2}$ and so on. Let $B^{\sharp}_{0}, B^{\sharp}_{1}, \dots$ be the resulting sequence of rescaled balls (so $B^{\sharp}_{0}=B_{0}$). If we reach stage $j$ then $\rho(q_{0})$ must be less than   $2^{-j}$ so,  since $\rho(q_{0})$ is  strictly positive, we must stop at some stage. Thus we have rescaled balls $B^{\sharp}_{0}, \dots B^{\sharp}_{J}$, all admitting $s$ charts. In the rescaled metric $d^{\sharp}$ on $B^{\sharp}_{J}$ the distance $d^{\sharp}(p_{J}, q_{0})$ is bigger than $1/10$ but less than $1/5+O(s)$ so we can suppose it is less than the quantity $\kappa_{1}- 2\kappa_{3}\kappa_{2}$ considered above. Thus we can apply Lemma 1 with initial condition $q=x_{0}= q_{0}$ The condition that $q_{0}$ lies in the image of $\Sigma$ implies that  $\rho(q_{0})$ is in roughly the same range.
In particular we can suppose that $\rho^{\sharp}(q_{0})$ is greater than $\kappa_{0}$ and less than $2\kappa_{2}$. Thus we can apply Lemma 1 with initial condition $q=x_{0}= q_{0}$, since obviously $q$ lies in $S_{q}$. Thus we can flow until $\rho^{\sharp}$ increases to $2\kappa_{2}$. Now pass to the ball $B^{\sharp}_{J-1}$, scaling by a factor of $2$. Our new initial condition has $\rho^{\sharp}=\kappa_{2}$ and we can again flow until $\rho^{\sharp}$ has increased to $2\kappa_{2}$. We continue this process through the sequence of balls $B^{\sharp}_{j}$ until in the original ball $B_{0}$ we flow to a point $x$ with $\rho(x) =2\kappa_{2}$ and $d(x,q)\leq \kappa_{3} \rho(x)$. This is not exactly what we asserted in Proposition 1, but it is clear that (by our assumption that $B_{0}$ lies in a larger ball which is close to the model) we can continue the process a little further until we hit $\chi_{s}(E)$.

\section{Completion of proof}

\subsection{Flat limit spaces: proof of Theorem 3}

Theorem 3 is not surprising and is very likely well-known to experts, but we did not find a precise statement in the literature. Notice that if  we cast our net wider to consider general length spaces which are the completions of flat Riemannian manifolds then  there is a huge amount of flexibility. For example, any $n$-dimensional simplicial complex with the property that all points lie in the closure of $n$-simplices can be given such a structure via an embedding in  Euclidean space. In particular the analogue of Theorem 3 would fail in this larger class.

Let $B$ be a flat limit ball as considered in Theorem 3. The first step is to see that we can suppose  that the s-chart $\chi_{s}:H_{s}\rightarrow B$ is actually an {\it isometry} (possibly after slightly shrinking the ball). Recall that if $V$ is any flat, connected, $n$-manifold  (not necessarily complete) there is a {\it developing map}: an open immersion of the universal cover of $V$ in $\bR^{n}$. This yields a developing homomorphism from $\pi_{1}(V)$ to the isometry group ${\rm Euc}_{n}$ of $\bR^{n}$. We apply this with $V=\chi_{s}(H_{s})$. Then $\pi_{1}(V)=\Gamma$ is finite and the homomorphisms from $\Gamma$ to ${\rm Euc}_{n}$ are rigid, up to conjugacy. This easily implies the statement.

The next observation is that, assuming Theorem 2 is known, there can be no codimension $2$ singular points---that is, no points with a tangent cone of the form $N_{\beta,n}$.
More generally no iterated tangent cone in $B$ can have this form.

 Write $H$ for the  the set of points $(\xi,\eta)\in \bR^{k}\times \bR^{n-k}$ such that $\vert \xi\vert \geq 1$.  
Let $\nu$ be a point in $S^{k-1}$, and let $r>0$. Define the subset
$U_{\nu,r}\subset H$ to be the points which can be joined to $(\nu,0)$ by a path {\it in $H$} of length at most $r$. This is a slightly complicated set but it is clear from the convexity of the complement of $H$ that the volume of $U_{\nu,r}$ strictly exceeds $\frac{\omega_{n}}{2} r^{n}$, where $\omega_{n}$ is the volume of the unit ball in $\bR^{n}$. (Note that $U_{\nu,r}$ is not quite the same as the intersection of $H$ with the ball $B(\nu,r)$, since $H$ is not convex.)  Fix $r$ small enough that $U_{\nu,r}$ and $g(U_{\nu,r})$ are disjoint for all $\nu$ and for all $g\in \Gamma$, different from the identity. Now fix $c>1$ such that
\begin{equation}   c^{n} \Vol(U_{\nu, r})> \frac{\omega_{n}}{2} (cr+c-1)^{n} \end{equation}

Thus $c$ depends only on $\Gamma$.  We can assume  that
\begin{equation}   \frac{1}{4s}> \left(1+ \frac{1}{c-1}\right)  \end{equation}  

For $0<\sigma<2s$ let $A_{\sigma}$ be the \lq\lq annulus $\{ \sigma <\vert \xi \vert< 2s\}$ in $\bR^{k}/\Gamma$. If $\eta\in \bR^{n-k}$ with $\vert \eta\vert$ not too large  we have an obvious isometric embedding
 $\iota_{\eta}:A_{s}\rightarrow B$ defined by the $s$-chart. Define $d(y)$ to be the infimum of the set of $\sigma$ such  that $\iota_{\eta}$ extends to an isometric immersion of $A_{\sigma}$ in the smooth part of $B$.
To simplify notation we suppose (as we obviously can by rescaling slightly) that $d(\eta)$ is actually defined for all $\eta$ with $\vert \eta \vert\leq 1 $.

\begin{lem}
Suppose that for all $\eta$ in some ball $D= \{ \eta: \vert \eta- \eta_{0}\vert <t\}$ we have $d(\eta)< \sigma $. Then the maps $\iota_{\eta}$ for $\eta$ in $D$ define   an isometric embedding of $D\times A_{\sigma}$ in $B$.
  \end{lem}
  This is certainly true when $\sigma=s$: the  map in question is just the restriction of the $s$-chart. We fix $\eta$ and let $\sigma_{\eta}$ be the infimum of the set of $\sigma$ such that the images  $\iota_{\eta}(A_{\sigma})$ are disjoint, as  $\eta$ ranges over $D$. If $\sigma_{\eta}>\sigma$ then there is a smooth point $q$ of $B$ which lies in the closure of $\iota_{\eta}(A_{\sigma_{\eta}})$ and of $\iota_{\eta'}(A_{\sigma})$ for some $\eta'\neq \eta$. (We allow the possibility that $\eta'$ is in the closure of $D$.) Elementary considerations involving the local geometry around $q$ show that this is impossible.(The essential point is the concavity of the boundaries of the sets in question.)

 With all this preparation we reach the central step in the proof.
\begin{prop}For $\eta_{0}$ with $\vert \eta_{0}\vert<1/2$ we have $d(\eta_{0})=0$\end{prop}

Suppose we have any point $\eta$ with $\vert \eta\vert<1$ and $d(\eta)=d>0$. The limit of the isometric embeddings $\iota_{\eta}(A_{\sigma})$ as $\sigma$ tends to $d$ defines a map from the closure $\overline{A_{d}}$ to $B$ which we still denote by $\iota_{\eta}$. Clearly there must be a singular point $q$ of $B$ in $\iota_{\eta}(\overline{A_{d}})$.
 Suppose, arguing for a contradiction,  that for all points $\eta'$ with $\vert \eta-\eta'\vert \leq d $ we have $d(\eta')\leq c d(\eta)$. For a suitable choice of $\nu$ there is an copy of $U_{\nu,r}$ in $B$ scaled by a factor $cd$ and contained in the ball of radius $d(cr+c-1)$ centred at $q$. Thus
(7) implies that
 
\begin{equation} \Vol(B(q,cdr))> \frac{\omega_{n}}{2} (cdr)^{n}. \end{equation}
But now we have
\begin{lem} Suppose $q$ is a singular point in an n-dimensional $L^{1}$-flat limit space $Z$. Then for any $d>0$ we have $\Vol(B(q,d)\leq \frac{\omega_{n}}{2} d^{n}$.
\end{lem}
Suppose first that there is a tangent cone at $q$ of the form $\bR^{n-p}\times C(X)$ with  $X$ smooth. Thus $X$ has the form $S^{p-1}/G$ where $G$ acts freely on the sphere, and the limit as $d$ tends to $0$ of $\frac{ \Vol(B(q,d))}{\omega_{n} d^{n}}$ is $1/\vert G\vert$. Then our result in this case follows from   generalised Bishop-Gromov monotonicity. 
Suppose next that there is a tangent cone at $q$ as above but with $X$ singular. Suppose however that there is a tangent cone of $C(X)$ of the  form $\bR^{n-p} \times \bR^{p}/G$, where $G$ acts freely. Then we first apply volume  monotonicity to see that $\Vol(X)\leq \omega_{n}/\vert G\vert$ and then argue as above. There is always some {\it iterated} tangent cone of the form $\bR^{n-p} \times \bR^{p}/G$, where $G$ acts freely, and we extend the argument in the obvious way. 

At this stage, combining Lemmas 2 and 3, we see that in fact if $d(\eta)>0$ there must be some point $\eta'$ with $\vert \eta'-\eta\vert \leq d(\eta)$ and $d(\eta')\geq c d(\eta)$. Start with a point $\eta_{0}$ with $\vert \eta_{0}\vert \leq 1/2$ and suppose $d(\eta_{0})=d>0$.
There is another point $\eta_{1}$ with $\vert \eta_{1}-\eta_{0}\vert\leq d$ and $d(\eta_{1})\geq c d$. Repeating the argument we get $\eta_{i}$ with
$\vert \eta_{i+1}-\eta_{i}\vert \leq d c^{i}$ and $d(\eta_{i})\geq c^{i} d$. We have to stop at some stage $\eta_{N}$,  when    we approach the boundary, so  $\vert \eta_{N}\vert \geq 1$. Then we have
\begin{equation}  d + c d + c^{2} d+\dots   + c^{N-1} d \geq \vert \eta_{0}-\eta_{N}\vert \geq 1/2. \end{equation}
On the other $d(\eta_{N-1})\leq s$ so $c^{N-1} d \leq s$. Combined with (10) this is a contradiction to our hypothesis (8). 

Given Proposition 6, we argue just as in Lemma 2 that the half-sized ball in $Z$ is isometric to a ball in $M_{\Gamma}$.

 Although we do not use this, Theorem 3 leads to a precise description of the local structure of $L^{1}$ flat limits.
\begin{prop}
An $L^{1}$ flat limit space is a Euclidean orbifold. Any point has a neighbourhood isometric to the quotient of the ball $B^{n}\subset \bR^{n}$ by a finite subgroup of $O(n)$.
\end{prop}

The holonomy of the flat metric defines a homomorphism $\pi_{1}(B_{\reg})\rightarrow
SO(n)$ and this in turn defines a covering of $B_{\reg}$. The basic point is that the metric completion of this covering is a smooth Riemannian manifold. It is easy to see that this implies Proposition 7. Now Theorem 3 has the following consequence. If there is an iterated tangent cone of $B$ of the form $M_{\Gamma,n}$ then there is a corresponding small ball in $B$ which is isometric to a ball centred at the origin in $M_{\Gamma,n}$. Now the assertion about the metric completion of the covering follows from an inductive argument on the codimension of the singular set which we  leave to the  interested reader to fill in. 

The work of Joyce \cite{kn:J} provides many examples of such Euclidean orbifolds, often with  intricate singular structures, which  arise as the limits of non-collapsed Einstein manifolds. 

There is a somewhat reciprocal relation between Theorem 3 and Proposition 7. In one direction, as we have sketched above, the first can be used to establish the second. Alternatively, if one knows Proposition 7 there are slightly simpler proofs of Theorem 3. 

 \subsection{Extension to bounded Ricci curvature}
 Recall that Theorem 1 deals with metrics of bounded Ricci curvature while in the previous sections we have restricted attention to the Einstein case.
Here we discuss the modifications required to remove this restriction. These are of a fairly  technical, but routine, nature. 

We recall that, in general, the Cromov-Hausdorff limit $W$ of a sequence of non-collapsed metrics of bounded Ricci curvature contains a dense open subset $W_{\reg}$  which is a manifold of class $C^{2,\alpha}$ (for any $\alpha\in (0,1)$) and the limiting structure on $W_{\reg}$ is a Riemannian metric of class $C^{1,\alpha}$. This is shown by working in harmonic co-ordinates on suitable small balls before passing to the limit. Thus the first modification is to change  the second item in the definition of an $s$-chart to say that the pull-back of $g$ differs in $C^{1,\alpha}$ norm from the flat metric by at most $s$. In fact all we need is $C^{1}$. Reviewing the proof we see that the only place where the regularity of the $s$-chart might be an issue occurs in Section 3.3, where we want to say that if the knot of a minimising lasso based at $x$ has a distance $O(s^{1/2})$ from $\varpi(x)$ then the tension vector differs from $\partial_{r}$ by $O(s^{1/2})$. But this is straightforward for metrics whose difference from the flat metric is small in $C^{1}$  since the  Christoffel symbols are then small in $C^{0}$. Thus solutions of the geodesic equation
$$   \ddot{x}_{i}= \Gamma^{i}_{jk} \dot{x}^{j}\dot{x}^{k}, $$
written in the $s$-chart, are close to Euclidean lines which easily yields what we need. 

The second issue arises in our discussion of $L^{1}$-flat limit balls.
We should now change the definition, replacing Einstein by bounded Ricci curvature. The Gromov-Hausdorff limit $Z$ still has an open dense subset $Z_{\reg}$ with a $C^{1,\alpha}$ metric, but it is not immediately clear that this  is, in suitable co-ordinates, a amooth flat Riemannian metric. To see this we go back to working in harmonic c-ordinates over small balls, before taking the limit. We recall that for $p>n/2$ there is a good theory of $L^{p}_{2}$ Riemannian metrics, with curvature in $L^{p}$. Elliptic estimates in harmonic co-ordinates mean that an $L^{\infty}$ (hence $L^{p}$) bound on the Ricci curvature give an $L^{p}_{2}$ bound on the metric on interior balls, and in particular an $L^{p}$ bound on the curvature. In our situation, the $L^{1}$ norm of the curvature tends to zero in the sequence, so for any $p'<p$ the $L^{p'}$ norm of the curvature tends to zero, by H\"older's inequality.
Fix $p'>n/2$, then we get an $L^{p'}_{2}$ limit, in harmonic co-ordinates, with curvature zero. Elliptic regularity implies that this is smooth. 
 Thus $Z_{\reg}$ has a smooth Euclidean structure and all our arguments carry over. 
\section{Codimension 2}

In this section we discuss extensions of  the argument we have given above  to establish  the codimension  2 situation;   Theorem 2.  We emphasise again that Theorem 2 is a known result and because of that  we will be content with a sketch. To simplify the discussion slightly we consider the Einstein case 

       The general set-up we considered in Section 2 and 3---the definition of an s-chart etc.---goes over in a obvious way. So now we write $H_{s}$ for the appropriate subset of $N_{\beta}= \bR^{2}_{\beta}\times \bR^{n-2}$. The new feature is that the analogue of Theorem 3 is false: the singularity is not rigid. Let $z_{1}, \dots z_{r}$ be points in $\bR^{2}$ such that we can find disjoint  \lq\lq wedges'' $W_{i}\subset \bR^{2}$ with vertices at the $z_{i}$ and with angles
$\gamma_{i}$. We get a flat singular space $Z$ by cutting out these wedges from the plane and gluing along the resulting edges. Then $Z$ has $r$ singular points with cone angles $2\pi(1-\gamma_{i})$ and outside a compact set is isometric to the standard cone with angle $\beta=1-\gamma$ where $\gamma=\sum \gamma_{i}$. Consider a family of such  wedges with fixed angles $\gamma_{i}$ and with vertices $z_{i}$ converging to the origin in $\bR^{2}$ and take the product with $\bR^{n-2}$.  We get  flat spaces which contain  isometric copies of $H_{s}$ for arbitrarily small $s$, but which are not isometric to
$N_{\beta,n}$. One can also construct more complicated examples where we have a unit ball in a singular flat space which contains an arbitrarily small deformation of $H_{s}$. This is possible because the fundamental group of $H_{s}$ is ${\bf Z}$ and the developing homomorphism 
from $ {\bf Z}$ to  ${\rm Euc}_{n}$  can be deformed. However it seems likely that the general picture will be similar to that considered above, with the singular set breaking up into a number of \lq\lq almost parallel'' copies with cone angles $\beta_{i}=1-\gamma_{i}$ where $\sum \gamma_{i}$ is as close as we please to $\gamma=1-\beta$. Thus to explain the argument simply we will  assume a rather weak version of this idea. 

\

{\bf Assumption} 

Fix $c>1/2$. Suppose there is an $s$-chart $\chi_{s}:H_{s}\rightarrow B$ in an $L^{1}$- flat limit ball. If $s$ is sufficiently small then either the half-sized ball $\frac{1}{2} B$ is isometric to a ball in $N_{\beta,n}$ or  there is a point
$q\in \frac{1}{2} B$ with a tangent cone $N_{\beta',n}$ where $(1-\beta')\leq c (1-\beta)$.

\

Now we argue as follows. Write $B^{n}$ for the unit ball in $\bR^{n}$. Recall that, according to Anderson and Colding, there is a $\delta_{0}$ such that if our unit ball $B(x,1)\subset X$ has Gromov-Hausdorff  distance less than $\delta_{0}$ from $B^{n}$ then we get a fixed bound on the curvature tensor in the half-sized ball. If the $L^{1}$ norm of the curvature tensor is small the same is true for all $L^{p}$ norms and then, by standard elliptic theory, for the $L^{\infty}$ norm. It follows that, given any $\delta_{1}$ we can find an $\epsilon(\delta_{1})$ so that if the $L^{1}$ norm of $\Riem$ is less than $\epsilon(\delta_{1})$ the Gromov-Hausdorff distance from the half sized ball $B(x,1/2)$ to $\frac{1}{2} B^{n}$ is less than $\delta_{1}$. Now suppose that $\beta$ is sufficiently close to $1$ that the distance from the unit ball in $B_{\beta,n}\subset M_{\beta,n}$ to $B^{n}$ is less than $\delta_{0}/2$. Set
$$\delta'= d_{GH}( \frac{1}{2} B_{\beta, n}, \frac{1}{2} B^{n}), $$
so $\delta'>0$, since $\beta<1$. Now if $\delta\leq \min(\delta_{0}/2, \delta'/2)$ it follows that we can {\it never} have the situation where $d_{GH}(B(x,1), B_{n,\beta})<\delta$ and the $L^{1}$ norm of the curvature is less than $\epsilon(\delta'/2)$.
This implies that  we cannot have a point in a flat limit space with tangent cone of the form $N_{n,\beta}$ for such values of $\beta$. Say this covers a range $1-\gamma_{0}\leq \beta<1$.

Now consider a value of $\beta$ with $1-\gamma_{0}/c<\beta<1$. Invoking the \lq\lq assumption'',  we see that a flat limit which admits an $s$-chart, for small enough $s$, is actually isometric to a ball in $N_{\beta,n}$ and our argument goes through.
 So we know that in fact we cannot have  a tangent cone in a flat limit space of the form $N_{n,\beta}$ for this larger range of $\beta$. Then repeat the argument until we cover the whole range $0<\beta<1$.

By making  more complicated arguments it seems that one can prove Theorem 2 using this approach, without using the \lq\lq assumption'' . This involves  another level of limits, taking $\delta\rightarrow 0$ and arguing with the resulting complete limit space.
For example, using  volume monotonicity and the \lq\lq volume cone implies metric cone'' theorem  one sees that if $\beta\geq 1/2$ then there are no higher codimension singularities in B, so certainly there must be {\it some} points with codimension 2 singularities. Similarly one sees that in proving the high codimension result, Theorem 1, it suffices to know Theorem 2 for the range $\beta\geq 1/2$. But the arguments become convoluted and, since the result is known, we leave the interested reader to fill in details. 

There is an unsatisfactory aspect in the statement of the \lq\lq assumption''. {\it A posteriori}---
given Theorem 2---none of these spaces can arise as $L^{1}$-flat limits. Further, if Anderson's conjecture in \cite{kn:A2} is true then none can arise as limits of metrics with bounded Ricci curvature. On the other hand they certainly do arise as limits of metrics with Ricci curvature bounded below Thus we feel that, properly formulated, there should be an interesting classification problem of a   deformations of the $N_{\beta,n}$ within a suitable class of flat singular spaces. 
 


\end{document}